\documentstyle[12pt]{amsart}
\newtheorem{lemma}{Lemma}

\newtheorem{theorem}{Theorem}

\newcommand{\R}{{\bf R}}
\newcommand{\C}{{\bf C}}
\newcommand{\rme}{{\rm e}}
\newcommand{\rmd}{{\rm d}}

\newcommand{\cH}{{\cal H}}

\newcommand{\alp}{\alpha}

\newcommand{\gam}{\gamma}
\newcommand{\kap}{\kappa}
\newcommand{\lam}{\lambda}
\newcommand{\del}{\delta}
\newcommand{\eps}{\varepsilon}

\newcommand{\lap}{{\Delta}}

\newcommand{\Dom}{{\rm Dom}}

\newcommand{\Spec}{{\rm Spec}}
\newcommand{\Ker}{{\rm Ker}}
\newcommand{\Ran}{{\rm Ran}}
\newcommand{\norm}{\Vert}
\renewcommand{\Re}{{\rm Re}\;}

\newcommand{\Proof}{\underbar{Proof}{\hskip 0.1in}}

\newcommand{\conv}{{\rm conv}}
\newcommand{\dist}{{\rm dist}}
\newcommand{\Schrodinger}{Schr\"odinger }
\newcommand{\lin}{{\rm lin}}
\newcommand{\la}{{\langle}}
\newcommand{\ra}{{\rangle}}
\newcommand{\pr}{\prime}
\newcommand{\EV}{eigenvalue \ }
\newcommand{\EVs}{eigenvalues \ }
\newcommand{\sqc}{\sqrt{c}}
\newcommand{\inn}{\int_{-\infty}^{\infty}}

\begin{document}
\title[Spectral Instability]{SPECTRAL INSTABILITY FOR SOME \\
SCHR\"ODINGER OPERATORS}
\author{A. Aslanyan and E. B. Davies}
\date{September 1998}
\address{Department of Mathematics, King's College,
Strand, London WC2R 2LS, UK}
\email{Aslanyan@@mth.kcl.ac.uk, E.Brian.Davies@@kcl.ac.uk}
\thanks{The authors thank the Engineering
and Physical Sciences Research Council
for support under grant No. GR/L75443}
\maketitle

\begin{abstract}
We define the concept of instability index of an isolated eigenvalue
of a non-self-adjoint operator, and prove some of its general
properties. We also describe a stable procedure for
computing this index for Schr\"o\-din\-ger operators in one
dimension, and apply it to the complex resonances of
a typical operator with a dilation analytic potential.
\end{abstract}


\section{Introduction}
In some earlier papers we showed that typical non-self-adjoint
\Schrodinger operators $H$ exhibit spectral instability in the following
sense. For any $\eps >0$ there exist many $\lam\in\C$ and
$f\in\Dom(H)$ such that
\[
\norm Hf-\lam f\norm_2 \leq \eps \norm f\norm_2
\]
even though $\lam$ is not near the spectrum of $H$. This
behaviour occurs for the harmonic oscillator with a nonreal coupling
constant as well as for many non-self-adjoint anharmonic oscillators.
There is a rapidly growing literature on pseudospectral theory,
which was invented to explore just such possibilities,
\cite{B1,B2,D1,D2,D3,Re,RedT,ReiT,TT,T1,T2}.

In this paper we return to the same type of
operator, but measure spectral instability by a method
which provides more precise information
about the instability of individual eigenvalues. We have
computed the so-called instability indices of the first
$100$ eigenvalues $\lam_n$ of the harmonic oscillator, and see that
they appear to increase exponentially with $n$. We have carried out a
similar but more limited exercise for the resonances of a typical \Schrodinger
operator with dilation analytic potential, and report our conclusions.

Our definition of the instability index of an isolated eigenvalue $\lam$
of $H$ of multiplicity $1$ involves the fact that the eigenfunction $f$
of $H$ associated with $\lam$ is different from the eigenfunction
$f^\ast$ of $H^\ast$ associated with its eigenvalue $\overline\lam$.
In Section 2 we show that the instability index
\[
\kap (\lam):=\frac{\norm f\norm_2 \norm f^\ast \norm_2}{|\la
f,f^\ast \ra |}
\]
of $\lam$ is equal to the norm of the spectral projection $P$ associated
with $\lam$, and also present a number of other theoretical
properties of the index.

If $H:=-\lap +V$ acts in $L^2(\R^N)$ where $V$ is a complex-valued
potential then $H^\ast=-\lap+\overline V$ and $f^\ast=\overline f$.
Hence for any isolated eigenvalue we have
\begin{equation}
\kap (\lam)=\frac{\int_{\R^{N}} |f|^2 }{|\int_{\R^{N}} f^2
|}.\label{index}
\end{equation}
Note that when  $\kap (\lam)$ is large the
eigenvalue is very unstable under small perturbations of the potential,
and hence also unstable because of rounding errors in the
computation.
The numerical task we face is to compute the eigenvalue and the
instability index accurately in situations in which the
denominator of (\ref{index}) is very small because the
complex-valued eigenfunction $f$ is oscillating rapidly.

Because the spectral instability develops so rapidly as $n$ increases,
we have had to take great care to use computational methods which
are reliable. Fortunately for our first problem there are independent
methods of checking the values which we have obtained. In the second
case we use the experience gained by the first problem, and have
checked the reliability of the conclusions under the variation of several
different parameters in the computational method. In Section~6 we
summarize the conclusions of our investigation.
\section{The Instability Index}

If $\lam$ is an isolated point of the spectrum of a closed
operator $A$ in a Hilbert space $\cH$, the spectral
projection $P$ associated with $\lam$ is defined by
\[
P\phi :=\frac{1}{2\pi i} \int_{\gam}(z -A)^{-1}\phi \rmd z
\]
where $\gam$ is any sufficiently small closed contour winding around $\lam$. The
assumption that this projection has rank $1$ is
stronger than the assumption that $\lam$ is an eigenvalue of
multiplicity $1$.

\begin{lemma}
If $f$ and $f^{\ast}$ are the normalised eigenvectors of $A$ and
$A^{\ast}$ associated with the eigenvalues $\lam$
and $\overline{\lam}$ respectively, and if $P$ has rank $1$,
then $\la f,f^{\ast}\ra \not= 0$ and
the instability index of $\lam$ is equal to $\norm P \norm$.
\end{lemma}

\Proof  We have
\begin{eqnarray*}
\Ker (P)& =&\left( \Ran(P^{\ast})\right)^{\perp}\\
&=& \{ g:\la g,f^{\ast}\ra =0\}.
\end{eqnarray*}
Since $f\notin \Ker (P)$ we see that $\la f,f^{\ast}\ra
\not= 0$. It is now easy to verify that $P$ is given by
\[
Ph:=\frac{ \la h,f^{\ast}\ra }{  \la f,f^{\ast}\ra}f
\]
and hence that
\[
\norm P\norm =\frac{1}{ |\la f,f^{\ast}\ra |} =\kap(\lam).
\]

\begin{theorem}
If $P$ has rank $1$ then
\[
\kap (\lam )=\sup \{ a_1(V):\norm V\norm \leq 1\}
\]
where $a_1(V)$ is defined to be the coefficient of $s$ in the expansion
of the eigenvalue of the perturbed operator $H+sV$:
\[
\lam(s)=\lam+a_1(V)s+a_2(V)s^2+\cdots
\]
\end{theorem}

\Proof  By standard arguments in perturbation theory
\cite{Ka} we have
\[
(H+sV)(f+sg+\ldots)=(\lam+s\mu+\ldots)(f+sg+\ldots)
\]
where the perturbed eigenvector is normalised by
\[
\la f+sg+\ldots,f^{\ast}\ra =\la f,f^{\ast}\ra.
\]
We deduce that $Hg+Vf=\lam g+\mu f$ and $\la g,f^{\ast}\ra
=0$. Therefore
\[
\la Hg,f^{\ast}\ra +\la Vf,f^{\ast}\ra= \mu \la f,f^{\ast}\ra
\]
and
\[
\mu=\frac{\la Vf,f^{\ast}\ra}{ \la f,f^{\ast}\ra}.
\]
The proof is completed by the observation that
\[
\sup\{ |\la Vf,f^{\ast}\ra|:\norm V\norm \leq 1\} =\norm f
\norm \,\norm f^{\ast}\norm.
\]

The instability index is also related to the notion of
pseudospectrum, which is a geometric way of looking at
resolvent norm properties. Namely if $\eps >0$ we put
\[
\Spec_{\eps}(A):=\{ z\in\C : \norm (z-A)^{-1}\norm
>\eps^{-1}\}.
\]
The sizes of these sets, which all contain the spectrum of
$A$, measure the spectral instability of $A$. One always has
\[
\{ z:\dist\{ z,\Spec (A)\} <\eps\} \subseteq \Spec_{\eps}(H)
\]
but the RHS  is often much larger than the LHS. The
following theorem states that if $\kap (\lam )$ is large then
the component of
$\Spec_{\eps}(A)$ containing $\lam$ is large in a related sense. Several
similar results can be obtained in the same manner.

\begin{theorem}
Suppose that the spectral projection $P$ associated with the isolated
eigenvalue $\lam$ of $A$ has rank $1$. Let
$\gam$ be a closed contour surrounding the
connected component of $\Spec_{\eps}(A)$ which contains
$\lam$ but does not intersect $\Spec_{\eps}(A)$. Then
\[
|\gam | \geq 2\pi\eps \kap (\lam)
\]
where $|\gam |$ is the length of $\gam$.
\end{theorem}

\Proof We have
\[
\kap(\lam)=\norm P\norm =\norm \frac{1}{2\pi i}
\int_{\gam}\frac{\rmd z}{z-A}\norm \leq
\frac{1}{2\pi\eps}\int_{\gam}|\rmd z|=
\frac{|\gam|}{2\pi\eps}.
\]
\section{The Computational Procedures}

\subsection{Finding Eigenvalues}

Let $H$ be the \Schrodinger operator
\begin{equation}
Hf(x):=-\frac{\rmd ^{2}f}{\rmd x^{2}} +V(x)f(x)
\label{sch}
\end{equation}
acting in $L^{2}(\R)$. Let us first outline briefly the method for
determining the eigenvalues of $H$. For any $z\notin [0,\infty)$
there exist two solutions $f_{\pm}$ of
\begin{equation}
Hf = zf \label{evp}
\end{equation}
which vanish as $x\to\pm \infty$ respectively.
We introduce {\em transfer functions} $g_{\pm}:=f_{\pm}^{\pr}/f_{\pm}$
satisfying  nonlinear first order differential equations
to be given below
and proper initial conditions at $\pm \infty$.
We wish to solve the Cauchy problem for $g_-(x)$ for
$x \geq X_-$ and for $g_+(x)$ for $x \leq X_+$
where $X_- < 0$, $X_+ > 0$, $|X_-|$ and $X_+$
are sufficiently large. Provided we
know initial conditions for $g_\pm$ at $X_\pm$
which correspond to $f_\pm$ vanishing at $\pm\infty$
respectively, the two Cauchy problems can be
solved by a standard numerical method.
The question how to transfer the so-called admissible boundary conditions
from singular points (which are $\pm \infty$ in our case) has been
extensively studied by A. A. Abramov and his collaborators
(a survey of their results can be found in \cite{AbrKon}).
Naturally, the behaviour of the
potential $V$ has to be taken into account: say, if $V$
is bounded and vanishes rapidly at infinity then
$g_\pm$ are asymptotically constant as $x\to\pm\infty$,
respectively. Later on in this section we shall consider
this and other cases applying the ideas of \cite{AbrKon}.

To locate the \EVs in terms of the transfer functions we
choose $a\in\R$ and consider
\[
F(z)=F(z;a):=g_{+}(a)-g_{-}(a).
\]
This function is meromorphic on $\C\backslash [0,\infty)$
with zeros at the eigenvalues of the operator $H$. The
zeros are independent of $a$ but $F$ may also have poles
which depend on $a$. They are at the eigenvalues of the
restrictions of $H$ to $L^{2}(a,\infty)$ and $L^{2}(-\infty,a)$ subject to
Dirichlet boundary conditions at $a$ in both cases.

To determine the zeros of $F$ numerically we use the argument
principle (cf. \cite{AbrYu} for a contour
integration procedure) to obtain their approximate positions followed by
some variant of Newton's method to obtain accurate values.
One has to be careful not to
choose a value of $a$ for which there is a pole close to
the zero of interest, so it is recommended that a few
different values of $a$ are investigated.

The numerical elaboration of the above ideas involves a
substantial amount of preliminary work. One may find the
approximate location of the eigenvalues and of the maxima
of the eigenfunctions by an independent method. For example
if one discretises a large enough interval in the real line then
one can find approximate eigenvalues by a standard matrix
eigenvalue routine; MATLAB is ideal for this purpose.
Another possibility is to use JWKB asymptotic formulae
which also enable one to find an interval $(X_{-},X_{+})$
outside which the relevant eigenfunction is negligible, and a
point $x_0$ at which the modulus of the eigenfunction takes its
maximum value (see the next subsection in this regard).

If the potential $V$ is an even function with
respect to reflection about the origin, then
every eigenfunction is either even or odd, and the problem
may be decomposed into two independent problems on
$[0,\infty)$, with Dirichlet or Neumann boundary conditions
at $0$. This is the case in our examples. From now on in this
subsection we concentrate on the symmetric problem
and take advantage of symmetry. Actual computational formulae
are given below for the case of the half-line. Note that the
theory of admissible boundary conditions
based upon asymptotic analysis of
solutions of the differential equation at $\pm\infty$
applies in a generic situation.

Auxiliary Cauchy problems to be solved numerically are as
follows. Let us first assume that $V(x) \sim cx^2$,
$x\to\infty$, where $c \in \C$ is constant --- this corresponds to
the harmonic oscillator problem
and its perturbations studied in the next section.
For a fixed value of $z$ we consider the solution $f_+$ of (\ref{evp})
vanishing as $x\to\infty$ and introduce a new transfer function $\alp_+$
satisfying
$$
\frac{1}{x} \alp_+' + \sqc \alp_+^2 + \frac{1}{x^2 \sqc} (\sqc \alp_+ - V
+ z) = 0, \qquad \alp_+(X) \sim -1 +
O\left(\frac{1}{X^2}\right), \qquad X\to\infty
$$
where $X$ is sufficiently large. According to \cite{AbrKon},
$$
f'_+(x)/f_+(x) =  \sqc x \alp_+(x)
$$
for such $x \geq 0$ that $\alp_+(x)$ exists.
Moreover, one can work out the
coefficients of the asymptotic expansion of $\alp_+$
for a particular potential $V(x)$. For $V(x) = cx^2$, say, we
replace the condition at infinity by
$$
\alp_+(X) = -1 + \frac{d}{X^2} +
\frac{\sqc d^2-d}{2 \sqc X^4}, \qquad d = \frac{z - \sqc}{2c},
$$
choosing $X$ so that the last term in the above formula
is negligible.
Along the lines of the mentioned paper we pose the
so-called admissible condition at infinity
for the considered potential. The above initial condition is
equivalent to the boundary condition $f_+ \to 0$ at infinity,
up to terms of order $O(\frac{1}{X^6})$ as $X \to\infty$.

Next, we are going to consider the operator $H$ with a potential
vanishing rapidly enough at infinity (see Section~5).
Following the same approach, we introduce
$$
g_+(x) = f'_+(x)/f_+(x), \qquad x \geq 0,
$$
and for this function obtain the singular Cauchy problem
$$
g_+^{\pr} + g_+^{2} - V + z = 0, \qquad g_+(x) \sim i \sqrt{z}, \qquad x\to\infty.
$$

Clearly, the initial conditions vary for different types of potential.
Still, for each particular choice of $V$ we are able to apply the
developed theory and set appropriate initial conditions for
transfer equations. After that has been done we integrate those equations
numerically from $X$ from right to left to some (fixed)
$a \geq 0$. The transfer
functions $\alp_+$, $g_+$ actually take their values in the
Riemann sphere, so we have to switch between them and their
inverses at certain values of $x$.
As soon as $|\alp_+|$ or $|g_+|$ becomes greater than a prescribed constant
we change to $\alp_+^{-1}$ or $g_+^{-1}$, respectively and from this point
on integrate analogous equations
starting with proper initial conditions for the inverse functions.
After a finite number of changes of this kind we reach the chosen
point $a$.

To complete the transfer procedure, consider the solution $f_0$ of
(\ref{evp}) satisfying $f'_0(0) = 0$.
We denote
$$
f_0'(x)/f_0(x) =  g_0(x)
$$
and solve
$$
g_0' + g_0^2 - V(x) + z = 0, \qquad g_0(0) = 0
$$
from left to right.
(Obvious changes should be made when considering an odd
solution $f_0$ satisfying a Dirichlet boundary condition at $x=0$.)
Again, we switch from $g_0$ to $g_0^{-1}$ if necessary and stop at the same
point $a$.

Remark that the described procedure is the simplest version of the
boundary condition transfer, or pivotal condensation method we have chosen for the second order
equation. There is an extensive literature
(cf. \cite{Abr, Bak, God} etc.)
on the transfer methods where much more advanced techniques are developed.
Although there are other possibilities, in this paper we
get satisfactory results implementing the above version.

Finally, we have $g_0$ (corresponding to either odd or even $f_0$)
and $g_+$ calculated at the same point $a$. If
$$
F(\lam;a) := g_+(a) - g_0(a) = 0
$$
then $\lam$ is an \EV of $H$. Thus, we evaluate $F(z;a)$
for certain values of $z$ as described and then
find the \EVs of interest as zeros of $F(z;a)$. As has already been mentioned,
we first use the argument principle (see \cite{AbrYu} for computational formulae)
to locate the \EVs and then apply an iterative Newton-like method to obtain more precise values.
After an \EV has been located up to the required accuracy, one can
compute the corresponding eigenfunction by recovering its values from
the transfer functions $g_+$, $g_0$.

The method proposed  has been implemented as a
universal Fortran~77 code including all the basic procedures
described above in this section.
Auxiliary Cauchy problems have been solved by a standard
routine based on the Runge--Kutta--Merson fourth order method.
In our computations we have used 32-bit and 64-bit arithmetic.

The question remains how to choose $a$ --- although the zeros of $F(z;a)$
do not depend on $a$, in practical computations the choice of $a$ does play a
significant role. In fact, we investigated different functions $F(z;a)$
for a wide range of $a$. One of the possible choices is
$a=0$. For problems with even potentials the zeros and the
poles of $F(z;a)$ interchange (see Figures~1a, 1b, 2a). We use
contour map plots to find initial guesses for \EVs when there is
no other a-priori information about their location.
Plotting contour maps with the use of Matlab~5.2 has also
helped us to avoid poles when looking for zeros. Thus,
in a generic situation we recommend that one
calculates the values of $F(z)$ (for which we have
used our Fortran code), then produces the
plots (Matlab graphics) and, finally, finds \EVs
accurately (a standard iterative Fortran procedure).

We regard our numerical results obtained
via the above method as reliable. In particular, this is confirmed
by several values of $a$ providing entirely different functions $F(z)$
whose zeros coincide. These results are to be reported below in the
following two sections.

\subsection{Calculating the Instability Index}

The second stage in the process is to compute the
instability index defined by (\ref{index}).
The obvious method, namely calculating the two
integrals after first determining the eigenfunction numerically,
is highly inaccurate if the instability index is large.
The reason is that the integrand in the denominator
is highly oscillatory, and the evaluation of
such integrals is problematical. The following method
is much superior in applications. Below we present a technique
suitable for an arbitrary (not necessarily symmetric) potential.

We introduce four functions as follows. For $x\in[a,X_{+}]$ we
define
\begin{eqnarray*}
        h_{+}(x) & := & f(x)^{-2}\int_{x}^{X_{+}} f(s)^{2}\rmd s \\
        k_{+}(x) & := & |f(x)|^{-2}\int_{x}^{X_{+}} |f(s)|^{2}\rmd s
\end{eqnarray*}
where $f$ is the eigenfunction associated with the
eigenvalue $\lam$.
Similarly for $x\in[X_{-},a]$ we define
\begin{eqnarray*}
        h_{-}(x) & := & f(x)^{-2}\int_{X_{-}}^{x} f(s)^{2}\rmd s \\
        k_{-}(x) & := & |f(x)|^{-2}\int_{X_{-}}^{x} |f(s)|^{2}\rmd s.
\end{eqnarray*}
It is obvious that
\[
\frac{k_{-}(a)+k_{+}(a)}{|h_{-}(a)+h_{+}(a)|}=
\frac{\int_{X_{-}}^{X_{+}}|f(x)|^{2}\rmd x}
{\left | \int_{X_{-}}^{X_{+}}f(x)^{2}\rmd x\right |}
\]
which converges exponentially rapidly
to $\kap$ as $X_{+}\to +\infty$ and  $X_{-}\to -\infty$.
The task is to find a procedure to evaluate the four
functions accurately. We consider only $h_{-}$, the others being similar.
It follows from its definition that $h_{-}(X_{-})=0$ and
that $h_{-}$ satisfies the differential equation
\begin{equation}
h_{-}(x)^{\pr}=1-2g_-(x)h_{-}(x).\label{hminus}
\end{equation}
This may be solved numerically, say, by a
Runge-Kutta method to determine $h_{-}(a)$.

It is important to be sure that the solutions of
(\ref{hminus}) and the other three equations are stable.
It suffices to note that $\Re g_+(X_+) < 0$ and $\Re g_-(X_-)
> 0$, which implies the stability of the solutions
$h_+$, $k_+$ and $h_-$, $k_-$ from right to left
and from left to right, respectively. This has been
confirmed by numerical tests. The same is true for the
transfer equations quoted in the previous subsection ---
the solution $\alp_+$, for instance, is known to be stable
from right to left which is essential for practical computations.

There is a potential problem in that if $f(b)=0$ for some
$b\in (X_{-},a)$ then $h_{-}$ is usually infinite at that point.
Generically one does not expect a complex-valued $C^{2}$ function of a
real variable to vanish anywhere, and we have not seen this
problem arise, but one needs to discuss how the method
should be adapted in the event of its occurrence. There are
two cases, which are distinguished numerically by whether
$h_{-}(x)\to 0$ as $x\to b$ or $|h_{-}(x)|\to\infty$ as
$x\to b$. Note that since $f$ is a non-zero solution of
(\ref{evp}) and $f(b)=0$ it follows that
$f^{\pr}(b)\not= 0$, so $|g(x)|\to \infty$ as $x\to b$.

\begin{lemma}
If $f(b)=0$ and $\int_{X_{-}}^{b}f(s)^{2}\rmd s=0$ then
$h_{-}(x)\to 0$ as $x\to b$ and $g(x)h_{-}(x)\to 1/3$ as
$x\to b$.
\end{lemma}

\Proof  Neglecting lower order terms we have
\begin{eqnarray*}
        h_{-}(x) & \sim &(x-b)/3   \\
        g(x) &\sim  &(x-b)^{-1}
\end{eqnarray*}
as $x\to b$. The results follow.

Thus, in this case we are still able to
integrate the same equation (\ref{hminus});
the point $b$  is, in fact, regular rather than
singular.

The more standard case is that in which $f(b)=0$ and
$\int_{X_{-}}^{b}f(s)^{2}\rmd s\not= 0$. Clearly $|h_{-}(x)|\to\infty$ as
$x\to b$.

\begin{lemma}
If we put $\tilde h_{-}(x):=h_{-}(x)^{-1}$ then $\tilde
h_{-}(x)\to 0$ and $\tilde h_{-}(x)g(x)\to 0$ as $x\to b$
and
\begin{equation}
\tilde h_{-}(x)^{\pr}=\tilde h_{-}(x)^2-2\tilde
h_{-}(x)g(x)   \label{htilde}
\end{equation}
for all $x$ near $b$.
\end{lemma}

\Proof Neglecting lower order terms we have
\begin{eqnarray*}
        \tilde h_{-}(x) & \sim &
        \frac{f^{\pr}(b)^{2}(x-b)^{2}}{\int_{X_{-}}^{b}f(s)^{2}\rmd s }\\
\tilde h_{-}(x)g(x)      & \sim &       \frac{f^{\pr}(b)^{2}(x-b)}{\int_{X_{-}}^{b}f(s)^{2}\rmd s }
\end{eqnarray*}
as $x\to b$. The verification that $\tilde h_{-}(x) $ satisfies the
differential equation (\ref{htilde}) is routine.

Thus, in the considered case we recommend to change to $\tilde{h}$
at $x=b$ and integrate (\ref{htilde}) instead of (\ref{hminus}).

Naturally, the stability of the procedure proposed in this subsection
depends heavily on $a$ (though the exact value of $\kap$ does not depend
on the norm of $f$). A proper choice of $a$ is very important and can
essentially influence the results. Choosing $a = {\rm argmax}
|f(x)|$ seems to be a reasonable way.

Compared to standard approaches the above mentioned technique
has two clear advantages.
First, we do not need to evaluate the fast oscillating
integrands  $f(x)^2$ and $|f(x)|^2$ themselves --- instead, we integrate
several auxiliary ODEs.  Secondly, this procedure is numerically
stable. In the cases which we
have examined the solutions $h_\pm$, $k_\pm$ change quite slowly and smoothly.

\subsection{Possible Difficulties}

If the instability index of an eigenvalue is very large then
it is clear from Theorem 2 that the eigenvalue is
intrinsically difficult to compute. One mechanism by which
this can occur in computations is that at the eigenvalues, for
which one knows that
$F(\lam)=0$, one also finds that $F^{\pr}(\lam)$ is very
small, so it is not possible to locate $\lam$ accurately. The
following theorem provides a link in one direction between these two
phenomena at a theoretical level.

We assume that
\[
H f(x):=-\frac{\rmd ^{2}f}{\rmd x^{2}} + V(x)f(x)
\]
on $L^{2}(\R)$, where $V(x)$ vanishes rapidly enough as $|x|\to\infty$.
Given $a\in\R$ and $z\in\C$ satisfying $\Re(i\sqrt{z})<0$, let $g_+(z,x)$
be the solution of
$$
g'(x) + g(x)^2 + z - V(x) = 0
$$
on $[a,\infty)$ subject to $g_+(z,x) \sim i\sqrt{z}$ as
$x\to\infty$. Let $g_-(z,x)$ be the solution of
the same equation on $(-\infty,a]$ subject to $g_-(z,x) \sim -i\sqrt{z}$ as
$x\to -\infty$. We put
$$
F(z) = g_+(z,a) - g_-(z,a)
$$
as usual so that $F(\lam)=0$ if and only if $\lam$ is an \EV of
$H$.

\begin{theorem}
Let $\kap(\lam)$ be the instability index at an
eigenvalue $\lam$, let $f$ be the associated eigenfunction
and assume that $g:=f^{\pr}/f$ is bounded on $\R$. Then
\[
\kap(\lam )|F^{\pr}(\lam) |\, \norm g\norm_{\infty} \geq 1.
\]
\end{theorem}

\Proof If $\eps >0$ is small enough there exists
\[
\mu=\lam+\frac{\eps}{F^{\pr} (\lam)} +O(\eps^{2})
\]
such that $F(\mu)=\eps$. Now put
\begin{eqnarray*}
        \tilde g_{-}(x) & := & g_{-}(\mu,x)+\eps/2  \\
        \tilde g_{+}(x) &:=  & g_{+}(\mu,x)-\eps/2
\end{eqnarray*}
for the appropriate values of $x$, so that
\[
\tilde g_{+}(a)-\tilde g_{-}(a)=0.
\]
We have
\begin{eqnarray*}
        \tilde g_{-}^{\pr}(x) & = & \tilde V(x)-\mu-\tilde g_{-}(x)^{2}  \\
        \tilde g_{+}^{\pr}(x) & = & \tilde V(x)-\mu-\tilde g_{+}(x)^{2}
\end{eqnarray*}
under the following conditions on $\tilde V$. If $x>a$ then
\begin{eqnarray*}
        \tilde V(x)-V(x) & = &\tilde g_{+}^{\pr}(x)+\mu +\tilde g_{+}(x)^{2}-V(x)   \\
         & = & g_{+}(\mu,x)^{\pr}+\mu +( g_{+}(\mu,x)-\eps/2)^{2}-V(x)  \\
         & = & -\eps  g_{+}(\mu,x)+\eps^{2}/4
\end{eqnarray*}
while if $x<a$ we must have
\begin{eqnarray*}
         \tilde V(x)-V(x)& = & \tilde g_{-}^{\pr}(x)+\mu +\tilde g_{-}(x)^{2}-V(x)  \\
         & = & g_{-}(\mu,x)^{\pr}+\mu +( g_{-}(\mu,x)+\eps/2)^{2}-V(x)  \\
         & = & \eps  g_{-}(\mu,x)+\eps^{2}/4
\end{eqnarray*}
Therefore
\[
\norm \tilde V-V\norm_{\infty} =\sup\{ \eps\norm
g_+(\mu,\cdot)\norm_{\infty}
+O(\eps^{2}),  \eps\norm g_-(\mu,\cdot)\norm_{\infty} +O(\eps^{2})\}.
\]
But $g_{\pm}(\mu,x)\to g(x)$ uniformly as $\mu\to\lam$ by the
assumptions of this section, so
\[
\norm \tilde V-V\norm_{\infty}=
\eps\norm g\norm_{\infty} +o(\eps).
\]
The statement of the theorem now follows from
the formula for the instability index
given in Theorem 2.

In the two examples considered below,
$F'(z)$ is very small for large values of $|z|$,
so it is impossible to determine its zeros.
This seems to be the main barrier to the
determination of large eigenvalues.

\section{The Harmonic Oscillator}

\subsection{Basic Facts}

Consider the operator $H$ defined by (\ref{sch}) with the
potential $V(x) = cx^2$, referred to as $H_o$ in the rest of the paper.
The eigenvalue problem for $H_o$ is called the
{\em harmonic oscillator} problem and is known to have infinitely many \EVs
$\lam_n^{(o)} = \sqc (2n+1)$, $n=0,1,\ldots$. The corresponding
eigenfunctions $f_n = C_n \rme^{-\sqc x^2/2} \phi_n(\sqrt[4]{c}\,x)$ where
$C_n$ are normalising constants, $\phi_n$ denote Hermite polynomials,
$\Re\sqc > 0$. These eigenfunctions are either even or odd:
$f_{2k}(x) = f_{2k}(-x)$, $f_{2k+1}(x) = - f_{2k+1}(-x)$, $k=0,1,\ldots$.
As proposed in Section~3, we consider $H_o$ on the half-line adding either
Neumann or Dirichlet boundary conditions at the origin.
We have used it as a sample problem to check the above method for finding
\EVs and eigenfunctions. Indeed, the results thus obtained are in very
good accordance with the theory; they confirm the reliability
of the method. It allowed us to calculate $\sim 100$ \EVs of $H_o$ up to the
accuracy $\del \in (10^{-10},10^{-4})$.

When implementing the method of Section~3 we found
that the accurate numerical determination of the \EV
$\lam_n$ for $n>100$ is not possible using double precision
(64-bit) arithmetic, particularly because $F'(z)$ can be very small near the points
where $F(z)=0$. We computed the instability indices for
the first 100 eigenvalues, using the JWKB approximation to the
eigenfunction $f_n$ associated with the \EV $\lam_n^{(o)}$
as described in \cite{D2, D3}. This approximation
suggests that $|f_n(x)|$ takes its maximum near $x=a_n$
and
$$
f_n(a_n+y) = \rme^{-i\eta_n y + O(y^2)}
$$
where the real constants $a_n$ and $\eta_n$ are
computed from
$$
\sqrt{c} (2n+1) = \eta_n^2 + c a_n^2.
$$
Having an appropriate value of $a_n$ is, of course,
helpful when we calculate $\kap_n$ by
means of the method given in Subsection 3.2.

\subsection{Perturbations of the Operator $H_o$}

Let us present some results concerned a perturbation of the
harmonic oscillator operator
$$
H_W := H_o + W(x).
$$
We have investigated various perturbations of the form
$W(x) = \eps \rme^{imx}$ for a range of
fairly small $\eps$. The reason is that
looking for the perturbation providing the most unstable
results, one has to choose $W(x)$ as follows.
From the perturbation theory formula cited in the proof of
Theorem~1 one can easily see that among  perturbations satisfying
$$
|W(x)| \leq \eps
$$
the function
$$
W_n(x)  = \eps \frac{\bar{f_n}(x)}{f_n(x)}
$$
provides the worst perturbation of the $n$-th \EV of $H_o$.
Indeed we then have
\begin{equation}
\lam_n = \lam_n^{(o)} + \eps \kap(\lam_n^{(o)})
 + o(\eps). \label{per}
\end{equation}

If we only take account of the first term of the
JWKB expansion for $f_n$, we obtain the perturbing
potential
$$
\tilde{W}_n(x) = \eps \rme^{2i\eta_n x}
$$
after removing an irrelevant phase factor. The
expectation that $\tilde{W}_n$ provides
a perturbation of the \EV almost as great as that due
to $W_n$ is tested below.
We tabulate below the absolute
values of the corrections to several \EVs $\lam_n$ of $H_W$
calculated numerically by means of the method described in Section~3.
Tables~1--3 contain the values of
$|\lam_n^{(o)} - \lam_n|$ for $c=\sqrt{i}$ and
$W(x) = \eps \rme^{imx}$. The figures related to $\eps = 0$
give the absolute errors of the computation of the
$n$-th \EV of $H_o$.

\pagebreak
Table~1. Values of $|\lam_n^{(o)} - \lam_n|, n = 9, \ 2\eta_n = 6.4133$

\vspace{.2in}

\begin{tabular}{|c|c|c|c|c|c|c|}
\hline
 $\eps \backslash m$
 & 1.0& 5.0& 6.0& 6.4133 & 7.0 & 10.0 \\  \hline
 0  &  $10^{-10}$ &   & & & & \\  \hline
 $10^{-6}$  & $7.9\cdot10^{-7}$ & $6.6\cdot10^{-6}$ &
 $7.7\cdot10^{-6}$ & $8.0\cdot10^{-6}$ &
 $7.8\cdot10^{-6}$ & $2.9\cdot10^{-7}$\\  \hline
 $10^{-5}$ & $8.6\cdot10^{-6}$ & $6.4\cdot10^{-5}$ &
 $7.9\cdot10^{-5}$ & $8.1\cdot10^{-5}$ &
 $7.7\cdot10^{-5}$ & $1.9\cdot10^{-6}$\\  \hline
 $10^{-4}$ & $8.7\cdot10^{-5}$ & $6.4\cdot10^{-4}$ &
 $8.0\cdot10^{-4}$ & $8.0\cdot10^{-4}$ &
 $7.7\cdot10^{-5}$ & $1.8\cdot10^{-5}$\\  \hline
 $10^{-3}$ & $8.7\cdot10^{-4}$ & $6.42\cdot10^{-3}$ &
 $7.97\cdot10^{-3}$ & $8.16\cdot10^{-3}$ &
 $7.66\cdot10^{-3}$ & $1.7\cdot10^{-4}$\\  \hline
\end{tabular}

\vspace{.2in}
Table~2. Values of $|\lam_n^{(o)} - \lam_n|, n = 19, \ 2\eta_n = 9.1884$

\vspace{.2in}

\begin{tabular}{|c|c|c|c|c|c|}
\hline
 $\eps \backslash m$ & 5.0& 9.0& 9.1884 & 10.0 & 20.0 \\  \hline
 0  &  $10^{-8}$  & & & & \\  \hline
 $10^{-8}$  & $1.0\cdot 10^{-7}$
 & $4.0\cdot 10^{-6}$ &$3.9\cdot 10^{-6}$ & $1.1\cdot 10^{-6}$ & $6.1\cdot 10^{-7}$\\  \hline
 $10^{-7}$ & $2.6\cdot 10^{-6}$
 & $3.25\cdot 10^{-5}$ & $3.22\cdot 10^{-5}$& $2.63\cdot 10^{-5}$ & $6.2\cdot 10^{-7}$\\  \hline
 $10^{-6}$ & $3.7\cdot 10^{-5}$
 & $3.25\cdot 10^{-4}$ &$3.20\cdot 10^{-4}$ & $2.77\cdot 10^{-4}$ & $6.1\cdot 10^{-7}$\\  \hline
 $10^{-5}$ & $4.0\cdot 10^{-4}$
 & $3.20\cdot 10^{-3}$ & $3.20\cdot 10^{-3}$ & $2.80\cdot 10^{-3}$ & $6.3\cdot 10^{-7}$\\  \hline
\end{tabular}

\vspace{0.2in}

Table~3. Values of $|\lam_n^{(o)} - \lam_n|, n = 29, \ 2\eta_n = 11.3014$

\vspace{0.2in}

\begin{tabular}{|c|c|c|c|c|c|}
\hline
 $\eps \backslash m$ & 10.0& 11.0& 11.3014& 12.0 & 20.0 \\  \hline
0  &  $10^{-6}$ &  & & & \\  \hline
$10^{-10}$ & $2\cdot 10^{-6}$ & $2\cdot 10^{-6}$
& $3\cdot 10^{-6}$& $3\cdot 10^{-6}$& $10^{-6}$\\  \hline
$10^{-9}$ & $1.2\cdot 10^{-5}$ & $1.5\cdot 10^{-5}$
& $1.4\cdot 10^{-5}$ & $1.4\cdot 10^{-5}$& $10^{-6}$\\  \hline
$10^{-8}$ & $1.14\cdot 10^{-4}$ & $1.48\cdot 10^{-4}$
&$1.46\cdot 10^{-4}$ & $1.34\cdot 10^{-4}$& $10^{-6}$\\  \hline
$10^{-7}$ & $1.141\cdot 10^{-3}$ & $1.454\cdot 10^{-3}$
& $1.441\cdot 10^{-3}$& $1.336\cdot 10^{-3}$& $10^{-6}$\\  \hline
$10^{-6}$ & 0.01142 & 0.01455
& 0.01453& 0.01334& $10^{-6}$\\  \hline
$10^{-5}$ & 0.11657 & 0.14927
& 0.14563& 0.13670& $10^{-6}$ \\  \hline
\end{tabular}
\vspace{0.2in}

First of all, the above results show that the
values of $|\lam_n^{(o)} - \lam_n|$ are approximately
proportional to $\eps$,
that is, confirm formula (\ref{per}) numerically.
In fact, for very small $\eps$ the method can only feel the first
order corrections to the \EVs of the harmonic oscillator within the
chosen accuracy as expected. Secondly, maximal perturbations of \EVs are
observed for $m \approx 2\eta_n$ which justifies the above
arguments.

\subsection{Another Approach. Instability Index}

We have also investigated the harmonic oscillator using the
quantum mechanical creation and annihilation operators
$A^{\ast}$ and $A$. This is not possible for generic
differential operators, but provides a method of testing the
general algorithms developed in the last section. In this language
\begin{eqnarray*}
        H_{o} & = & P^{2}+cQ^{2}  \\
         & = & -(A^{\ast}-A)^{2}/2 +c(A^{\ast}+A)^{2}/2 \\
         &=& (c-1)A^{\ast 2}/2 +(c+1)A^{\ast}A +(c-1)A^{ 2}/2
         +(c+1)/2.
\end{eqnarray*}

If $\{\phi_{n}\}_{n=0}^{\infty}$ is the orthonormal basis
of Hermite functions in $\cH := L^{2}(\R)$, then $A\phi_{n}=\sqrt{n}\phi_{n-1}$
and $A^{\ast}\phi_{n}=\sqrt{n+1}\phi_{n+1}$ for all $n$,
and we may represent $H_{o}$ by means of the infinite matrix
\[
H_{o,m,n}:=\left\{ \begin{array}{ll}
a_{m}& \mbox{ if $m=n$}\\
b_{m}&\mbox{if $n=m+2$}\\
b_{n}&\mbox{if $m=n+2$}\\
0&\mbox{otherwise}
\end{array}
\right.
\]
with respect to this basis, where $m,n=0,1,2,\ldots$
and
\begin{eqnarray*}
        a_{m} & := & (c+1)(m+1/2)  \\
        b_{m} & := & (c-1)\{(m+1)(m+2)\}^{1/2}/2.
\end{eqnarray*}

The even and odd subspaces $\cH_{0}$ and $\cH_{1}$ with
respect to reflection about $0$
are invariant under $H_{o}$, and these subspaces may also
be characterised by
\begin{eqnarray*}
        \cH_{0} & = & \lin\{\phi_{2n}:n=0,1,\ldots\}  \\
        \cH_{1} & = & \lin\{\phi_{2n+1}:n=0,1,\ldots\} .
\end{eqnarray*}
Restricting the matrix to either of these subspaces renders
it tri-diagonal, so numerical computations are particularly
easy and accurate. We compute the instability index of an
eigenvalue $\lam_{r}=\sqrt{c}(2r+1)$ for $r=0,1,\ldots$ by evaluating
$$
\kap_r:=\frac{\sum_{n=0}^{N-1}|f(n)|^{2}}{|\sum_{n=0}^{N-1}f(n)^{2}|}
$$
where $f$ is the eigenvector associated with $\lam_{r}$,
obtained by solving the obvious recurrence relation
starting from $n=0$. Note that $\lam_{r}$ is taken to be an exact
eigenvalue of the infinite matrix, not an
eigenvalue of the truncated $N\times N$ matrix. For a particular
eigenvalue $\lam_{r}$, $N$ must be large enough for the coefficients
$f_{n}$ with $n>N$ to be insignificant, but not so large that
the recurrence relation becomes unstable. For $r>50$ it is
not possible to satisfy both of these conditions
simultaneously using
standard double precision $32$-bit arithmetic, and we used
the high precision arithmetic of Maple V.4.

The delicacy of the computations is indicated by the
evaluation of $\kap_{100}\sim 2.5594\times 10^{16}$ for
$c=\sqrt{i}$. For $N=200$ this required us to use
the command `$Digits=30$' in
Maple V.4, but putting $N=500$, we only obtained the same
result for `$Digits=110$' or greater. The instability in the
solution of the recurrence relation is evidently more
important than the contributions of the terms of the series
in the range $200<n<500$. The following results (see Table~4)
were all obtained with $N=200$ and
`$Digits=100$', and appear to be reliable.

The instability indices tabulated below have been
obtained in two independent ways.
The methods developed in this and the previous sections
turned out to provide very close results for the first
40 eigenvalues. This can be seen from Table~4 where
$\kap_n^{(1)}$ is related to the method of this section,
and $\kap_n^{(2)}$ to that of Section~3.
The figures obtained for $n\geq 40$ are
clearly different for the two methods although they
are qualitatively of the same order.

\vspace{.2in}

Table 4. Instability indices of $H_o$, $c=\sqrt{i}$

\vspace{.2in}

\begin{tabular}{|c|c|c|c|c|c|c|}
\hline
 $n$  &  0 & 10 & 20 & 30 & 40 & 50  \\ \hline
$\kap_n^{(1)}$  &  1.0404 & 14.2777 & 563.2146 &
  2.5789$\cdot 10^4$ & 1.2625$\cdot 10^6$ &
  6.3627$\cdot 10^7$\\ \hline
$\kap_n^{(2)}$  &  1.0404 & 14.2777 & 563.2146 &
  2.5789$\cdot 10^4$ & 1.2625$\cdot 10^6$ &
  6.3649$\cdot 10^7$\\ \hline
$n$ & 60 & 70 & 80 & 90 & 100 & \\ \hline
$\kap_n^{(1)}$ &  3.2734$\cdot 10^9$&  1.7081$\cdot 10^{11}$
 & 9.0059$\cdot 10^{12}$ & 4.7860$\cdot 10^{14}$
 &2.5594$\cdot 10^{16}$ & \\ \hline
$\kap_n^{(2)}$ &  3.2922$\cdot 10^9$&  1.7110$\cdot 10^{11}$
 & 8.9063$\cdot 10^{12}$ & 4.0052$\cdot 10^{14}$
 &1.9261$\cdot 10^{16}$ & \\ \hline
\end{tabular}

\vspace{.2in}

The growth of the instability index corresponds to the values
of $|\lam_n^{(o)} - \lam_n|$ increasing with $n$ (see also
Tables~1--3). Results to be cited below provide another
numerical evidence of this fact. In Table~5 the values of
$|\lam_n^{(o)} - \lam_n|$, $c=\sqrt{i}$, corresponding to
the perturbing potentials $\tilde{W}_n$ are given.
Comparing Table~5 to Table~4 we conclude
that $\mu_n \approx \kap_n$
which indicates reasonably good agreement of our numerical results and
perturbation theory.

\vspace{.2in}

Table~5. Values of $|\lam_n^{(o)} - \lam_n|$, $\tilde{W}_n(x) = \eps \rme^{2i\eta_n
x}$, $c=\sqrt{i}$

\vspace{.2in}

\begin{tabular}{|c|c|c|c|c|}
\hline
 $\eps \backslash n$& 30& 40 & 50 & 60\\  \hline
$10^{-6}$ & 0.021542& 1.19860& & \\  \hline
$10^{-7}$ & 0.002155& 0.10747& & \\  \hline
$10^{-8}$ & 0.000216& 0.01056& 0.54950& \\  \hline
$10^{-9}$ & $2.2 \cdot 10^{-5}$& 0.00105& 0.05518& 2.4921 \\  \hline
$10^{-10}$& $3 \cdot 10^{-6}$& 0.00010& 0.00587 & 0.2587 \\  \hline
$10^{-11}$ & & $10^{-5}$& 0.00059& 0.0279\\ \hline
$10^{-12}$ & & & $6 \cdot 10^{-5}$& 0.0028\\ \hline
0 & $10^{-6}$& $5 \cdot 10^{-6}$& $10^{-5}$& $10^{-4}$\\ \hline
\end{tabular}

\vspace{.2in}

Analysing the rate of divergence of the instability index
of $H_o$ in Table~4, one can notice that it grows exponentially:
$\kap_n \sim \rme^{0.4 n}$ for the studied range of $n$.

The eigenfunctions for the harmonic operator with nonreal coupling constant do
not form an unconditional basis, \cite{D2}. If they
formed a conditional basis the projections $P_n$
associated with the eigenvalues $\lam_n=\sqrt{c}(2n+1)$ as in
Lemma 1 would be uniformly
bounded in norm by a standard argument, \cite{GK}. However, we have
obtained strong numerical evidence that the norms increase exponentially
with $n$. We therefore make the conjecture that for nonreal coupling constant
the eigenfunctions of the harmonic oscillator do not form a
conditional basis.

More precisely let $N >0$ and let $P_{N}$ be the
spectral projection of $H_o$ associated with the first $N$ complex
eigenvalues $\lam_{n}^{(o)}$ where they are ordered in
increasing absolute values.  Explicitly
\[
P_{N}f:=\sum_{n=0}^{N-1}\frac{\la f,f_{n}^{\ast}\ra}{\la
f_{n},f_{n}^{\ast}\ra}f_{n}.
\]
If $\norm P_{N}f-f\norm\to 0$ as $N\to+\infty$ for all
$f\in\cH$ then the
uniform boundedness theorem implies that there exists a
constant $C$ such that $\norm P_{N}\norm \leq C$ for all
$N$. From the inequality
\[
\norm P_{N}-P_{N-1}\norm \leq 2C
\]
we are then able to deduce that the instability index
$\kappa_{n}$ is a bounded function of $n$. This conflicts with
the numerical evidence that these indices increase
exponentially with $n$. We have attempted to confirm the
exponential increase by using the JWKB approximations to the
eigenfunctions constructed in \cite{D2, D3}, but the eigenfunctions
oscillate so rapidly for high eigenvalues that the JWKB
approximations were not useful. While we have not proved the
exponential increase of $\kappa_{n}$ the corresponding result for the
pseudospectrum (resolvent norms) has been proved in
\cite{D3} not just for the harmonic oscillator but for a
wide range of anharmonic oscillators.

\section{Complex Resonances}

\subsection{Definitions}

Let $H$ be the \Schrodinger operator
\[
\hat{H}f(x):=-\frac{\rmd ^{2}f}{\rmd x^{2}} + V(x)f(x)
\]
acting in $L^{2}([0,\infty))$ subject to Dirichlet or Neumann
boundary conditions at $x=0$,
where the potential $V$ is bounded and
vanishes at infinity. For any positive constant $c$ we
define
\begin{eqnarray}
\hat{H}_{c}f(x) & := & \left( D_{c}\hat{H} D_{c^{-1}}f\right)
(x)\nonumber\\
&=&-c^{-2}\frac{\rmd ^{2}f}{\rmd x^{2}} +V(cx)f(x) \label{Hc}
\end{eqnarray}
where $D_{c}$ is the unitary dilation operator
\[
D_{c}f(x):=\sqc f(cx).
\]
We observe that $\hat{H}_{c}$ is unitarily equivalent to $\hat{H}$.
If $V$ is an entire function on $\C$ then
the formula (\ref{Hc}) defines a family of
non-self-adjoint operators
parametrised by $c\in\C$, $c\not= 0$. Under suitable conditions
the eigenvalues of these operators
are known to be independent of $c$, and are called
{\em resonances} of $\hat{H}$; see \cite{CFKS, HS} for expositions
of the theory of dilation analytic resonances. Since the
operators $\hat{H}_{c}$ are unitarily equivalent for values of $c$
with the same argument, we only consider $c$ of the form
$c:=\rme^{i\theta /2}$ where $0<\theta <\pi/2$.

We investigate the particular case of the operator
\[
H_0f(x):=-\frac{\rmd ^{2}f}{\rmd x^{2}} +x^{2}\rme^{-x^{2}/b^{2}}f(x)
\]
where $b>0$ is to be fixed. If one imposes a Dirichlet
boundary condition at $x=0$, this operator determines the
evolution in the zero angular momentum sector of a
three-dimensional quantum particle trapped by a rotationally
invariant barrier, where the
particle may tunnel through the barrier and escape to
infinity. Because the potential is non-negative and vanishes
rapidly at infinity, $H_0$ has absolutely continuous
spectrum $[0,\infty)$ and no eigenvalues.
A direct calculation shows that
\[
H_{\theta}f(x):= \hat{H}_{\rme^{i\theta/2}}f(x) = -\rme^{-i\theta}\frac{\rmd ^{2}f}{\rmd x^{2}}
+\rme^{i\theta}x^{2}\rme^{-\rme^{i\theta}x^{2}/b^{2}}f(x).
\]
We consider $H_\theta$ subject to either Dirichlet or Neumann
boundary conditions at $x=0$.
The potential of this operator vanishes rapidly as $x\to\infty$
provided $0<\theta<\pi/2$. Under this condition $H_{\theta}$
has essential spectrum $\rme^{-i\theta}\times [0,\infty)$ and also
some isolated eigenvalues in the sector $\{ z:-\theta <\arg z
\leq 0 \}$, these being independent of $\theta$.

For large values of $b$ (we take $b=100$) the potential of
$H_\theta$ is similar to that of the complex harmonic oscillator, and
the \EVs of $H_{\theta}$ are close to the values
$\{2n+1: \ n=0,1,\ldots\}$. For smaller values of $b$
there are several resonances very close to the positive real
axis, but at a certain point they turn sharply away into the
lower half plane.

\subsection{Location of Resonances}

The reason for there being resonances very close to the
real axis is as follows. Let us consider the operator
$H_\theta$ as a perturbation of the harmonic oscillator
operator $H_o$. In our notation we now have
$$
H_\theta = H_o + W(t, \nu), \qquad
W(t, \nu) = t^2(\rme^{-t^2/b^2} - 1)
$$
where we put $t = \rme^{i\theta/2}x$ and
$\nu = 1/b^2$. Regarding $\nu$ as a small parameter
we expand
$$
W(t, \nu)
= \sum_{k=1}^{\infty} W_k(t) \nu^k
= \sum_{k=1}^{\infty} \frac{(-1)^k t^{2(k+1)}}{k!} \nu^k.
$$
Again, for an arbitrary $n$,
following the standard perturbation theory approach, we
expand the $n$-th \EV of $H_\theta$ as
\begin{equation}
\lam_n(\nu) \sim \lam_n^{(o)} + \sum_{k=1}^{\infty} \mu_k \nu^k \label{as}
\end{equation}
which is a non-convergent asymptotic expansion,
and calculate
\begin{equation}
\mu_1 = \frac{\inn W_1(t)f_n^2 \rmd t}{\inn f_n^2 \rmd t} =
- C_n^2 \inn t^4 \rme^{-t^2} \phi_n^2(t) \rmd t. \label{fir}
\end{equation}
Here we follow the notations of Section 4: $f_n$ are the
eigenfunctions of $H_o$ and $\phi_n$ are Hermite
polynomials.

Formula (\ref{fir}) implies
that the first order correction $\mu_1$ is real and does
not depend on $\theta$. The same is true for all $\mu_k$.
Indeed, it is easily seen that the parameter $\theta$ enters
the problem in a specific way.
If one passes to the new variable $t$ and
proceeds with calculation of higher order corrections,
all the relations thus obtained do not contain any complex values
except for $t$ as an integration variable. Thus, one only deals with
integrals of the form $\inn p(t) \rmd t$ which do not depend on
$\theta$ and, therefore, are real.

Using the  creation--annihilation technique based on the
corresponding decomposition of the operator $H_o$
(see the previous section), we calculate the first order
correction for the
$n$-th \EV implicitly. Thus, formula  (\ref{fir}) becomes
\begin{equation}
\mu_1 = \mu_1(n) = - \frac{3}{4} (2n^2+2n+1),
\qquad n=0,1,\ldots. \label{first}
\end{equation}

Following the same numerical
procedure (see Section~3) we compute some of
the \EVs of $H_\theta$.
These results can be found in Subsection~5.4.

\subsection{Numerical Range and Complex Resonances}

The resonances must turn away from the real axis as their
absolute value increases, because of the fact
that a resonance $z$ is an eigenvalue of $H_{\theta}$. This
implies that
\[
z\in\bigcap_{0<\theta<\pi/2}N(\theta)
\]
where $N(\theta)$ is the numerical range of the operator $H_{\theta}$.
The numerical range of $H_{\theta}$ is defined by
\begin{eqnarray*}
        N_{\theta} & := & \{ \la H_{\theta} f,f\ra:\norm f \norm =1\}  \\
         & = & \{ \int\{ V_{\theta}(x)|f(x)|^{2}
         +\rme^{-i\theta}|f^{\pr}(x)|^{2}\} \rmd x : \norm f \norm  =1\} \\
         & \subseteq & \{ \int  V_{\theta}(x)|f(x)|^{2} \rmd x:\norm f \norm  =1\}
         +\rme^{-i\theta}[0,\infty)\\
         & \subseteq  & \conv \{ V_{\theta}(x):x\in \R\} +
         \rme^{-i\theta}[0,\infty),
\end{eqnarray*}
where
\[
V_{\theta}(x):= \rme^{i\theta}x^{2}\rme^{-\rme^{i\theta}x^{2}/b^{2}}.
\]
For small positive $\theta$ the set $N(\theta)$
crosses the real axis near $x\sim 0.461 b^{2}$.

Indeed, if we represent $V_\theta(x) = X_\theta(x) + Y_\theta(x)$
and denote the point where $N(\theta)$ meets the real axis by $(B_\theta,0)$
then a simple calculation gives us
$$
B_\theta = X_\theta + Y_\theta \cot\theta = b^2
\max\{x^2\rme^{-x^2}(2-x^2): \ x \in \R\}
+ O(\theta) = 0.461 b^2 + O(\theta).
$$
Therefore the imaginary parts of any resonances must start
decreasing before their real parts reach this value. This is
in good accordance with the numerical data quoted in the next
subsection.

\subsection{Numerical Results}

Lower \EVs of the operator $H_\theta$ lying close to
the real axis for  different values of $\nu$
are given in Tables~6,~7. The computed
\EVs proved not to depend on $\theta$, so
our numerical results are in agreement with the theoretical
arguments of Section~5.1. The fact that for a range of $\theta$
lower \EVs coincide up to a high accuracy shows
the stability of our method as a whole.

Remark that the results to be
reported below are consistent with formulae (\ref{fir}) and
(\ref{first}); they confirm, in particular, that $\mu_1 < 0$.
On the other hand, these results illustrate the fact that
series (\ref{as}) is asymptotic rather than convergent. This
only implies that the imaginary parts of the resonances
have to be very small within the regime for which the
asymptotic expansion provides useful information.

\vspace{.2in}

Table~6. Resonances of $H_0$, $n=0, \ 1$

\vspace{.2in}
\begin{tabular}{|c|c|c|}
\hline
 $\nu$  &  $\lam_0$ & $\lam_1$ \\ \hline
  0. & 1. &  3. \\ \hline
 $10^{-4}$ & 0.999925 &  2.999677 \\ \hline
 $4\cdot10^{-4}$ & 0.999700 & 2.998502 \\ \hline
 $10^{-3}$ & 0.999251 & 2.996253 \\ \hline
 $10^{-2}$ & 0.992475 & 2.962115 \\ \hline
 0.04 & 0.969405 &  2.824312 \\ \hline
 0.1 & $0.920295 -2.\cdot 10^{-6} i$  & $2.560861 - 0.003347 i$ \\ \hline
 0.2 & $0.822647 - 0.005282 i$ & $2.028250 -0.249944 i$ \\ \hline
 0.25 & $0.768023 - 0.019417 i$ & $1.850388 - 0.425748 i$ \\ \hline
\end{tabular}

\vspace{.2in}

Table~7. Resonances of $H_0$, $\nu=10^{-4}$

\vspace{.2in}

\begin{tabular}{|c|c|c|c|}
\hline
$n$ & $\lam_n$ & $n$ & $\lam_n$\\ \hline
0 & 0.999925 &  18 & 36.948571 \\ \hline
1 &  2.999677 & 20 & 40.936844 \\ \hline
2 &  4.999025 & 22 & 44.923925  \\ \hline
3 &  6.998125 & 24 & 48.909797  \\ \hline
4 &  8.996924 & 26 & 52.894462 \\ \hline
5 &  10.995877 &28 & 56.877922  \\ \hline
6 &  12.993623 &30 & 60.860175  \\ \hline
7 &  14.991344 &40 & 80.753321  \\ \hline
8 &  16.989119 & 50 & 100.616298 \\ \hline
9 &  18.986242 & 60 & 120.449097 \\ \hline
10 &  20.983418 &  70 & 140.237430 \\ \hline
12 & 24.976502  & 80 & 160.000434 \\ \hline
14 & 28.968399  & 90 & 179.874965  \\ \hline
16 & 32.959086  & 100 & 199.664121 \\ \hline
\end{tabular}

\vspace{.2in}

Along with Table~7 we present some
plots (see Figures 1a--1c). They include contour maps
of the function $F(z;a)$ defined in Section~3 whose zeros are the
\EVs we are looking for. One can see that the zeros and
the poles of $F$ interchange (we have plotted $F$
for Neumann boundary condition at $x=0$,
i.e., its zeros are the even eigenvalues, while the poles
correspond to the odd ones). For $|z| \leq 200$ we have
discovered 100 \EVs all being real up to the accuracy
$\del = 10^{-6}$. Note that for
different values of $a$ (an intermediate matching point) we obviously get quite
different functions $F(z;a)$ (compare Figure~1a to 1c) while their zeros remain the
same.

Given a certain number $n$ we watch $\lam_n$
changing as $\nu$ increases and compare this \EV with
its first order approximation (see (\ref{as})).

\vspace{.5in}

Table~8. Eigenvalues of $H_\theta$, $n = 10, \ 30$

\vspace{.2in}

\begin{tabular}{|c|c|c|c|c|}
\hline
 $\nu$  &  $\lam_{10}$ & $\lam_{10}^{(o)} + \mu_1(10)\nu$
 &  $\lam_{30}$ & $\lam_{30}^{(o)} + \mu_1(30)\nu$  \\ \hline
$10^{-4}$ & 20.9834  &  20.9834 & 60.8602 & 60.8604 \\ \hline
$10^{-3}$ & 20.8332  &  20.8343 & 59.5781 & 59.6043 \\ \hline
$2\cdot 10^{-3}$ & 20.6643  &  20.6685 & 58.0976 & 58.2085 \\ \hline
$4\cdot 10^{-3}$ & 20.3195  &  20.3370 & 54.9090 & 55.4170 \\ \hline
$6\cdot 10^{-3}$ & 19.9644  &  20.0005 & 51.2492 & 52.6255 \\ \hline
$8\cdot 10^{-3}$ & 19.5531  &  19.6740 & $46.2597 - 0.1492 i$& 49.8340 \\ \hline
$9\cdot 10^{-3}$ & 19.3747  &  19.5083 & $43.3597 - 1.5496 i$& 48.4383 \\ \hline
$10^{-2}$ & 19.2045  &  19.3425 & $41.2601 - 3.2488 i$ & 47.0425 \\ \hline
\end{tabular}

\vspace{0.2in}

It is seen from Table~8 that only for
a narrow range of $\nu$
do the perturbation theory formulae (\ref{as})
approximate actual \EVs (compare to Table~6).
As $\nu$ increases a typical \EV deviates gradually
from the value given by (\ref{as}) and
at some stage its imaginary part becomes substantial.

The values of the instability indices of resonances
depend on $\theta$, even though the positions of the
resonances do not. We have observed that the
indices are in fact monotonically increasing functions of
$\theta$. While this is not surprising we have no proof of
the fact. We have also observed that the instability indices
are increasing functions of $b$, provided one follows the
`same' resonance as $b$ increases. The instability
indices $\kap_n, \ n=0, \ 10, \ 20$, computed
for a wide range of $\theta$ and $b$ are given
in the following three tables.

\pagebreak

Table~9.1. Instability indices $\kap_n, \ n=0$

\vspace{.2in}

\begin{tabular}{|c|c|c|c|c|c|c|}
\hline
 $\nu$ $\ \backslash \theta$ &
 $\pi/30$ &$\pi/8$ & $\pi/6$ & $\pi/4$ & $\pi/3$ & $\pi/2.5$
 \\ \hline
 0. & 1.002750 & 1.040381 & 1.074570& 1.189207 &
 1.414214 & 1.798908 \\ \hline
$10^{-4}$ & 1.002750 & 1.040378 & 1.074563 & 1.189185 &
 1.414134 & 1.798588 \\ \hline
$10^{-2}$ & 1.002729 & 1.040039 & 1.073886 & 1.186951 &
 1.406329 & 1.769120 \\ \hline
\end{tabular}

\vspace{.2in}

Table~9.2. Instability indices $\kap_n, \ n=10$

\vspace{.2in}

\begin{tabular}{|c|c|c|c|c|c|}
\hline
 $\nu$ $\ \backslash \theta$ &
 $\pi/3000$ & $\pi/300$ & $\pi/30$ & $\pi/20$ & $\pi/10$ \\ \hline
 0. & 1.000028 & 1.0233 & 1.3299 & 1.8249 & 6.6784 \\ \hline
 $10^{-4}$ & 1.000027 & 1.0209 & 1.3294 & 1.8243 & 6.6728 \\ \hline
 $10^{-2}$ & 1.000031 & 1.0028 & 1.3046 & 1.7562 &5.9928 \\ \hline
 $\nu$ $\ \backslash \theta$ &
 $\pi/8$ & $\pi/6$ & $\pi/5$ & $\pi/4$ & $\pi/3$ \\ \hline
 0. & 14.2777 &  57.4539 & 195.9499 & 1565.2614
 & $1.3645551 \cdot 10^5$ \\ \hline
 $10^{-4}$ & 14.2836 &
 57.3505 & 195.4619 & 1558.9429 & $1.3507237 \cdot 10^5$ \\ \hline
 $10^{-2}$ & 12.3306 & 45.7594 & 143.8155 & 965.0957 &
 $4.6122845 \cdot 10^4$ \\ \hline
\end{tabular}

\vspace{.2in}

Table~9.3. Instability indices $\kap_n, \ n=20$

\vspace{.2in}

\begin{tabular}{|c|c|c|c|c|c|c|}
\hline
 $\nu$ $\ \backslash \theta$ &
 $\pi/20$ & $\pi/10$ & $\pi/6$ & $\pi/5$ & $\pi/4$ & $\pi/3$ \\ \hline
 0. & 5.9275 & 113.5766& 9850.7214 & $1.1753 \cdot 10^5$&
 $7.4538\cdot 10^6$ & $3.6609\cdot 10^{10}$ \\ \hline
 $10^{-4}$ & 5.9190 & 113.1898 & 9782.2812 & $1.1641 \cdot 10^5$&
 $7.3412\cdot 10^6$ &
 $3.6119\cdot 10^{10}$ \\ \hline
 $10^{-3}$ & 5.6055 & 109.2017 & 9128.5324 & $1.0610 \cdot 10^5$&
 $6.3605\cdot 10^6$ &
 $3.4647\cdot 10^{10}$ \\ \hline
 $10^{-2}$ & 4.7349 & 75.0929 & 4811.8348 & $4.5459 \cdot 10^4$&
 $1.6894\cdot 10^6$ &
 $1.7152\cdot 10^9$ \\ \hline
\end{tabular}

\vspace{.2in}

Finally, let us cite some results obtained for $b=10$.
In this case we have found numerically several
resonances,
which are real up to the chosen accuracy $\delta =
10^{-4}$, and a series of complex ones. As is seen,
starting from about $n=20$ their imaginary part rapidly
increases in absolute value. In fact, for different
values of $\theta$ the number of resonances with negative
imaginary parts varies. The resonances and the relevant
instability indices are tabulated below.
In Table~10 we cite the \EVs $\lam_n$ of
$H_\theta$ along with the corresponding instability
indices $\kap_n$ calculated for $\theta=\pi/4$ and $\theta=\pi/16$.

The data of Table~10 is illustrated by
the plot of $F(z)$ (Figure 2a) and its
contour maps (Figures 2b,~2c). Remark that the
largest instability indices for $\theta=\pi/4$ and $\theta=\pi/16$
correspond to the 26-th and the 24-th \EV respectively
(here we have concentrated on even eigenvalues; odd \EVs
behave similarly). Figures~2b,~2c also show that the most
unstable \EVs are related to $n=24$ and $n=26$. They
appear to be the first \EVs with negative imaginary parts
--- as one can see, the following \EVs go to the complex
plane quite abruptly. We do not have any theoretical
explanation of this fact except for the remark on the
boundedness of $\kap_n$ made in the end of Section~4.
Anyway, the contour maps and the values of $\kap_n$
agree very well and imply the same --- the maximum of
$\kap_n$ is obtained for the `critical' range of the
spectral parameter where \EVs start moving away from
the real axis.

Note that though for $\theta=\pi/16$ and $\theta=\pi/4$
the contour map plots are quite similar, this only
means that the first 28 \EVs coincide for the two
operators. As we know, there are no \EVs of $H_\theta$
below the line $\rme^{-i\theta}\times [0,\infty)$.
The spots indicating the zeros of the function
$F(z)$ which are
beyond the range $\{ z:-\theta<\arg z\leq 0 \}$
correspond to solutions of $H_\theta f = z f$ growing
at infinity rather than decaying. Thus, in the considered example
we should only regard the first 28 zeros as the \EVs of $H_\theta$,
$\theta=\pi/16$. They coincide with those obtained for
$\theta=\pi/4$ as we expected. We believe that our results
are reliable because of their stability under the
variation of several parameters involved in the
problem.

\vspace{.2in}

Table~10. Values of $\lam_n$ and $\kappa_n$ for $\nu=10^{-2}$

\vspace{.2in}

\begin{tabular}{|c|c|c|c|}
\hline
$n$ & $\lam_n$ & $\kap_n, \ \theta=\pi/4$ & $\kap_n, \ \theta=\pi/16$ \\ \hline
0 & 0.9925 & 1.1870 & 1.0097 \\ \hline
2 &  4.9009 & 2.8983 & 1.0684 \\ \hline
4 &  8.6836 & 11.3609& 1.2100\\ \hline
6 &  12.3350 & 49.4772& 1.4462\\ \hline
8 &  15.8488 & 219.4180 & 1.7974\\ \hline
10 &  19.2174 & 960.5058 & 2.2918 \\ \hline
12 & 22.4312 & 4075.82 & 2.9652\\ \hline
14 &25.4782 & $1.6515\cdot 10^4$ & 3.8576 \\ \hline
16 &28.3422 & $6.2860\cdot 10^4$ & 5.0033\\ \hline
18 &31.0004 & $2.1989\cdot 10^5$ & 6.4039\\ \hline
20 &$33.7512 - 0.0003 i$ & $1.3978\cdot 10^6$ & 9.5706\\ \hline
22 & $35.5098 - 0.0014 i$ & $1.5650\cdot 10^6$ & 10.0018\\ \hline
24 & $37.0693 - 0.1593 i$ & $2.4535\cdot 10^6$ & 10.2337\\ \hline
26 & $38.7468 - 1.0004 i$ & $2.8963\cdot 10^6$ & 8.8755\\ \hline
28 & $39.8045 - 2.0367 i$ & $2.5627\cdot 10^6$ & 7.0743\\ \hline
30 & $41.2601 - 3.2488 i$ & $1.8339\cdot 10^6$ &  \\ \hline
32 & $42.6021 - 4.7565 i$ & $1.3337\cdot 10^6$ &  \\ \hline
34 & $45.6102 - 8.0230 i$ & $4.2730\cdot 10^5$ &  \\ \hline
36 & $47.0034 - 9.8515 i$ & $2.3134\cdot 10^5$ &  \\ \hline
\end{tabular}

\vspace{.2in}

\section{Conclusions}

The instability index of an eigenvalue of a non-self-adjoint ordinary
differential operator was
defined in Section 2, where we investigated its theoretical
properties. We have described a known general numerical procedure for
computing the eigenvalues of the differential operator and have introduced
a new and numerically stable procedure for computing the instability
indices.

In order to test this procedure, we have carried out extensive
computations for the harmonic oscillator with a complex coefficient. The
eigenvalues $\lam_n$ of this operator are given by an exact formula, and we found
close agreement between the formula and our numerical results for $n\leq
100$. We have also computed the instability index by two independent
methods, the first being the general procedure mentioned above. The second
uses a special numerical technique only available for the harmonic operator,
but capable of yielding extreme accuracy if implemented in Maple with
high precision arithmetic. The instability indices of the first
$40$ eigenvalues obtained by the two methods were found to be in close
agreement; see Table~4.

The discrepancies between the two methods are partly explained by the very high
values of the instability indices of the eigenvalues $\lam_n$ for $n\geq 40$.
This phenomenon was first observed for the harmonic oscillator in \cite{D2,D3}
where we approached the phenomenon via pseudospectral theory. Our current
approach has the advantage that it provides a quantitative measure of the
instability of individual eigenvalues under small perturbations of the
potential. We have carried out numerical experiments and confirmed that
the size of the effects predicted matches what we have observed for a
particular perturbation. In Section~4 we have conjectured on the
basis of the numerical results that the eigenfunctions of the
complex harmonic oscillator do not form a conditional basis.

We have also investigated the complex resonances of a typical self-adjoint
operator by means of the standard technique of dilation analyticity. This
identifies the resonances of the original operator with eigenvalues of any one
of a family of associated non-self-adjoint operators indexed by an angle.
The eigenvalues of these operators are independent of the angle, but the
instability indices depend upon its value.

We have discovered that for a certain operator, as is seen from Table~10,
the first $20$ eigenvalues have very small imaginary parts,
which is explained by the fact that
there exists a (non-convergent) asymptotic expansion which has real coefficients of all
orders. For higher eigenvalues the imaginary parts of the eigenvalues
increase rapidly in absolute value. We have computed the instability
indices of these eigenvalues for typical angles and discovered that
they increased rapidly with the modulus of the eigenvalue,
reaching a maximum value near the region where the imaginary part starts to increase
(see Tables~9.1--9.3,~10). No theoretical explanation of this phenomenon exists.

The very large size of the instability indices in both examples
indicates
that the computation of large eigenvalues of non-self-adjoint differential
operators is likely to be intrinsically intractable in many other cases of
a similar type. The same applies to the computation of large
resonances of self-adjoint differential operators.
The effect of rounding errors or of small
perturbations of the operator may be to change the
computed eigenvalues drastically. This discovery casts some doubt
on the significance of theoretical investigations of the
asymptotic distributions of resonances or of any computations of such
eigenvalues for all except self-adjoint operators. Our experience, and
that of others who work within the pseudospectral approach, has been that
the extreme instability of large eigenvalues is the norm rather than
a possibility which occurs only in pathological cases.

\vskip 0.5in
{\bf Acknowledgments} The authors are grateful to Prof.
A. A. Abramov, who first suggested this approach to
the investigation of instabilities, for useful
discussions of different aspects of the paper.
\par


\eject
\begin{figure}[h]
\vspace{1cm}
\begin{picture}(200,200)(0,0)
\includegraphics{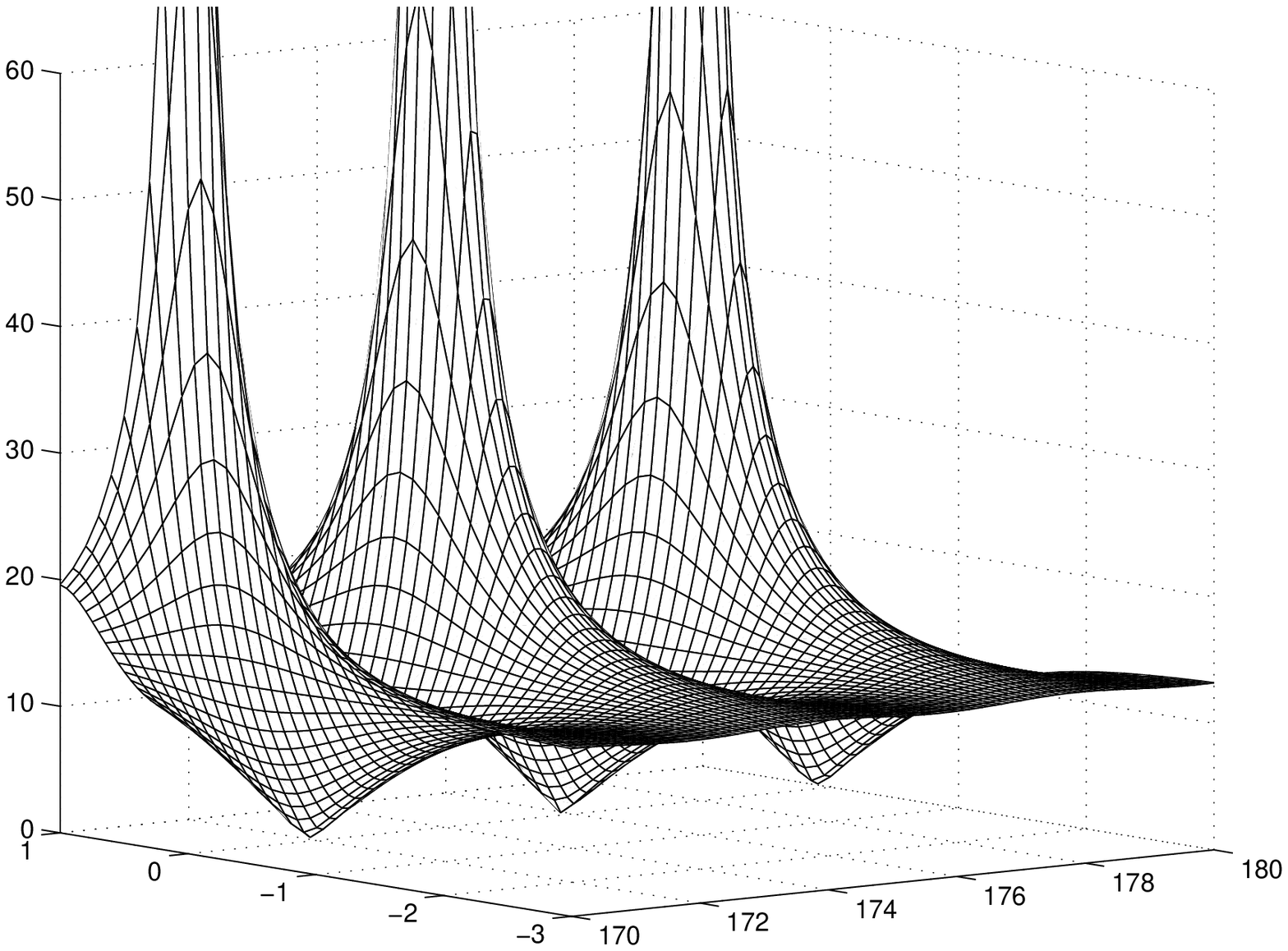}
\end{picture}
\end{figure}

\begin{center}
Figure 1a.  Function $F(z;a)$, $a=0$, $b=100$
\end{center}

\begin{figure}[h]
\begin{picture}(200,200)(0,0)
\includegraphics{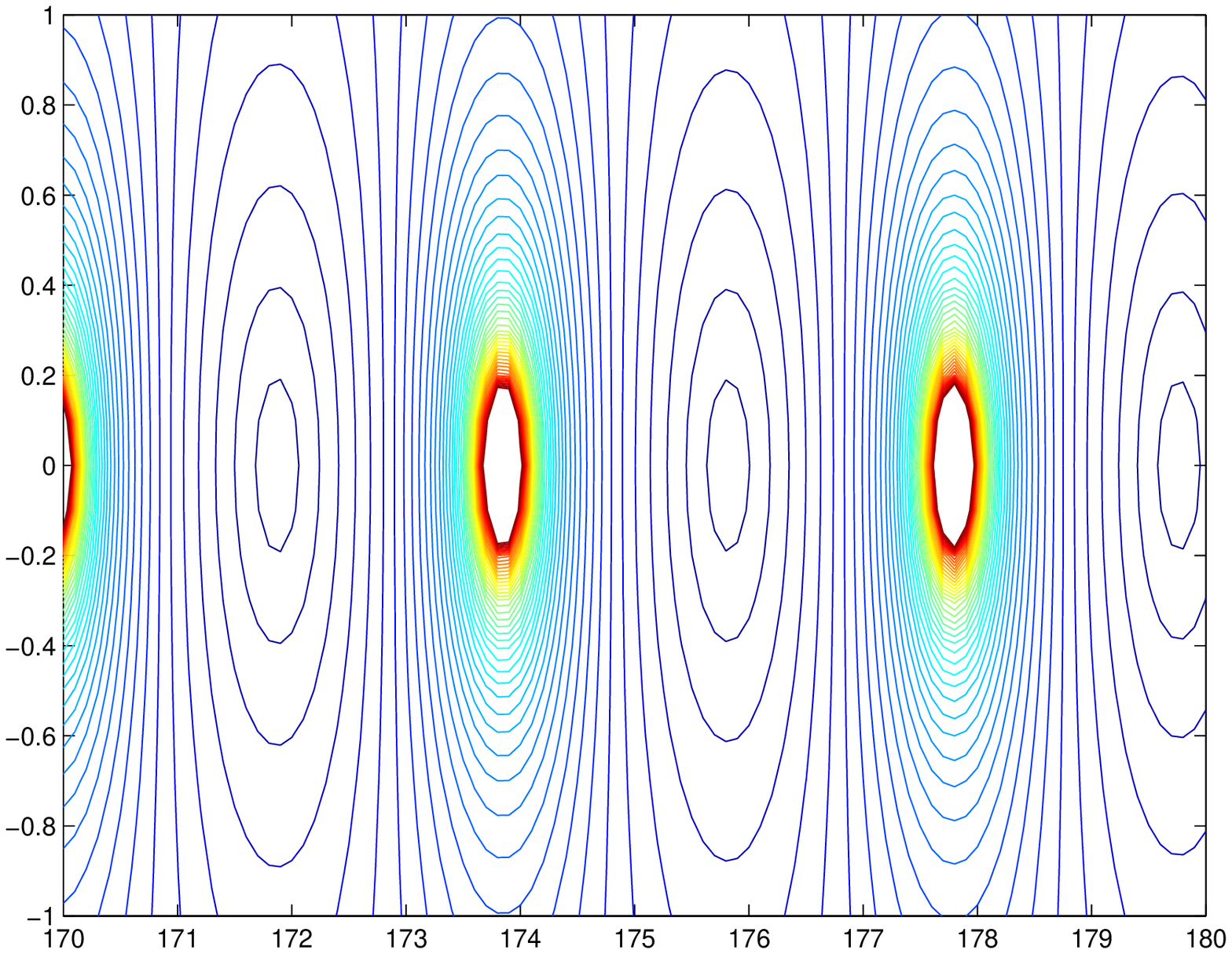}
\end{picture}
\end{figure}

\begin{center}
Figure 1b. Contour map of $F(z;a)$
\end{center}

\eject
\begin{figure}[h]
\vspace{1cm}
\begin{picture}(400,400)(0,0)
\includegraphics{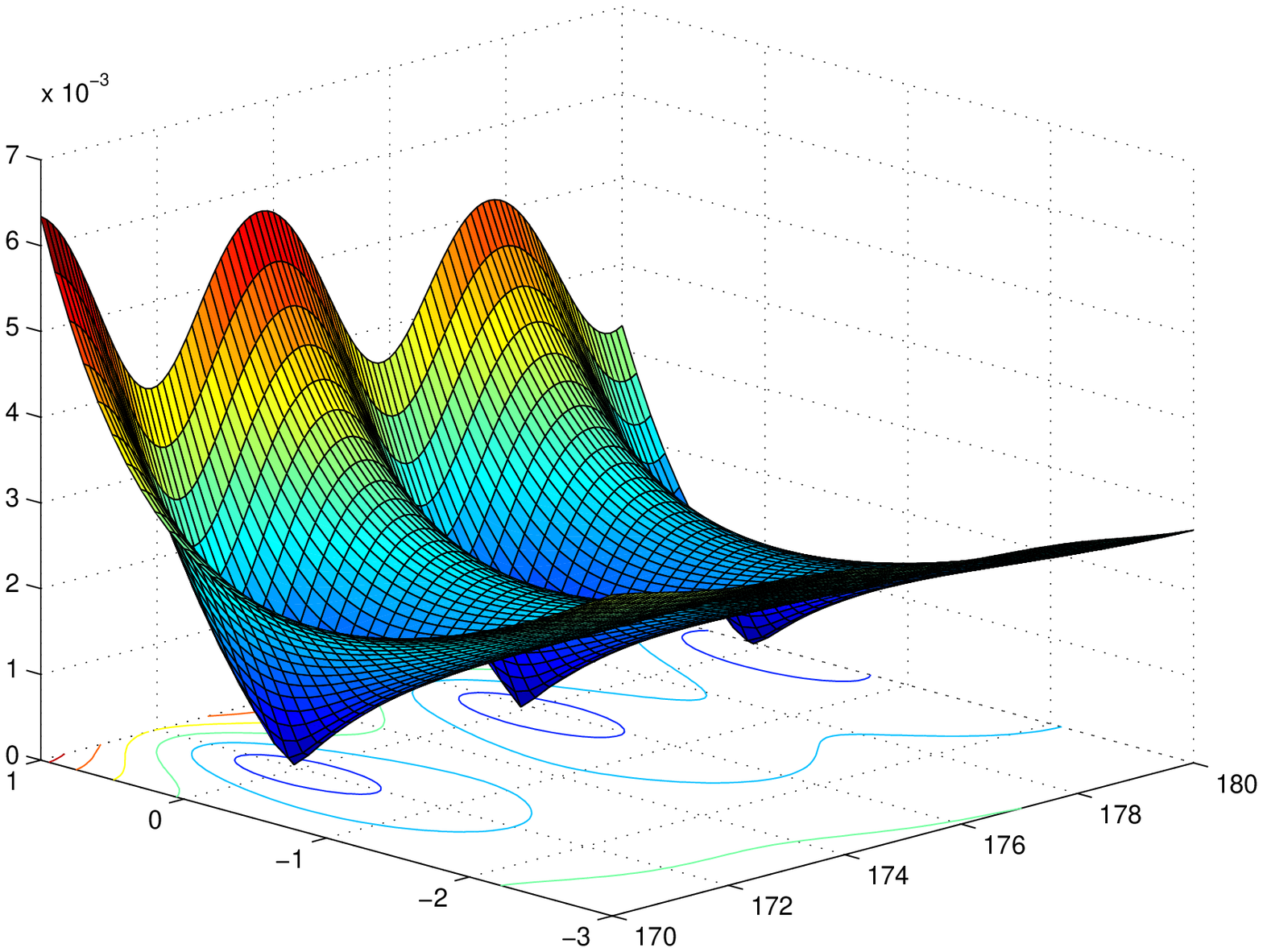}
\end{picture}
\end{figure}

\vspace{-2in}
\begin{center}
Figure 1c.  Plot and contour map of $F(z;a)$, $a=4$, $b=100$
\end{center}

\eject

\begin{figure}[h]
\vspace{1cm}
\begin{picture}(400,400)(0,0)
\includegraphics{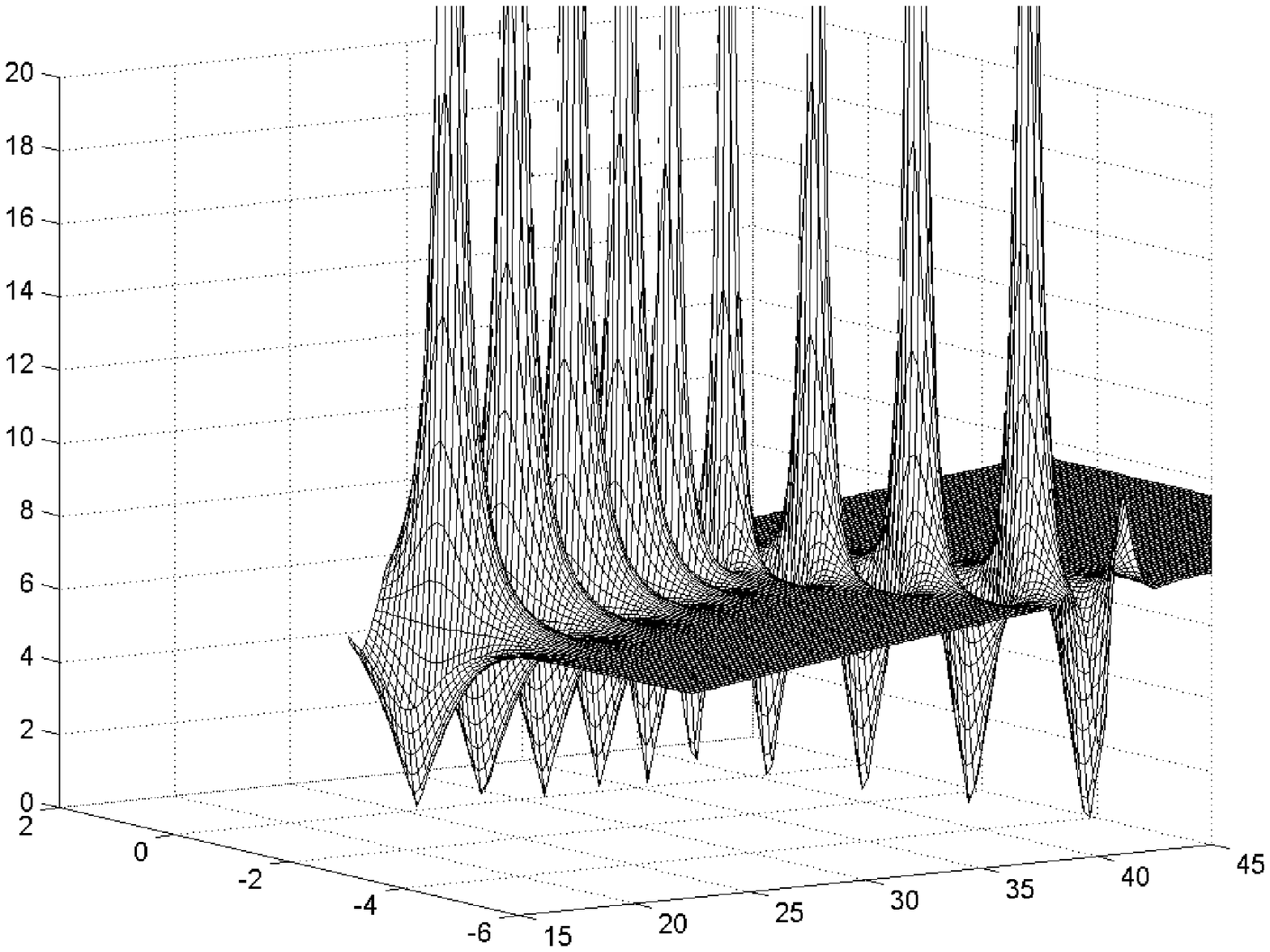}
\end{picture}
\end{figure}

\vspace{-2in}
\begin{center}
Figure 2a. Function $F(z;0)$, $b=10$, $\theta=\pi/4$
\end{center}

\eject

\begin{figure}[h]
\vspace{-1cm}
\begin{picture}(200,200)(0,0)
\includegraphics{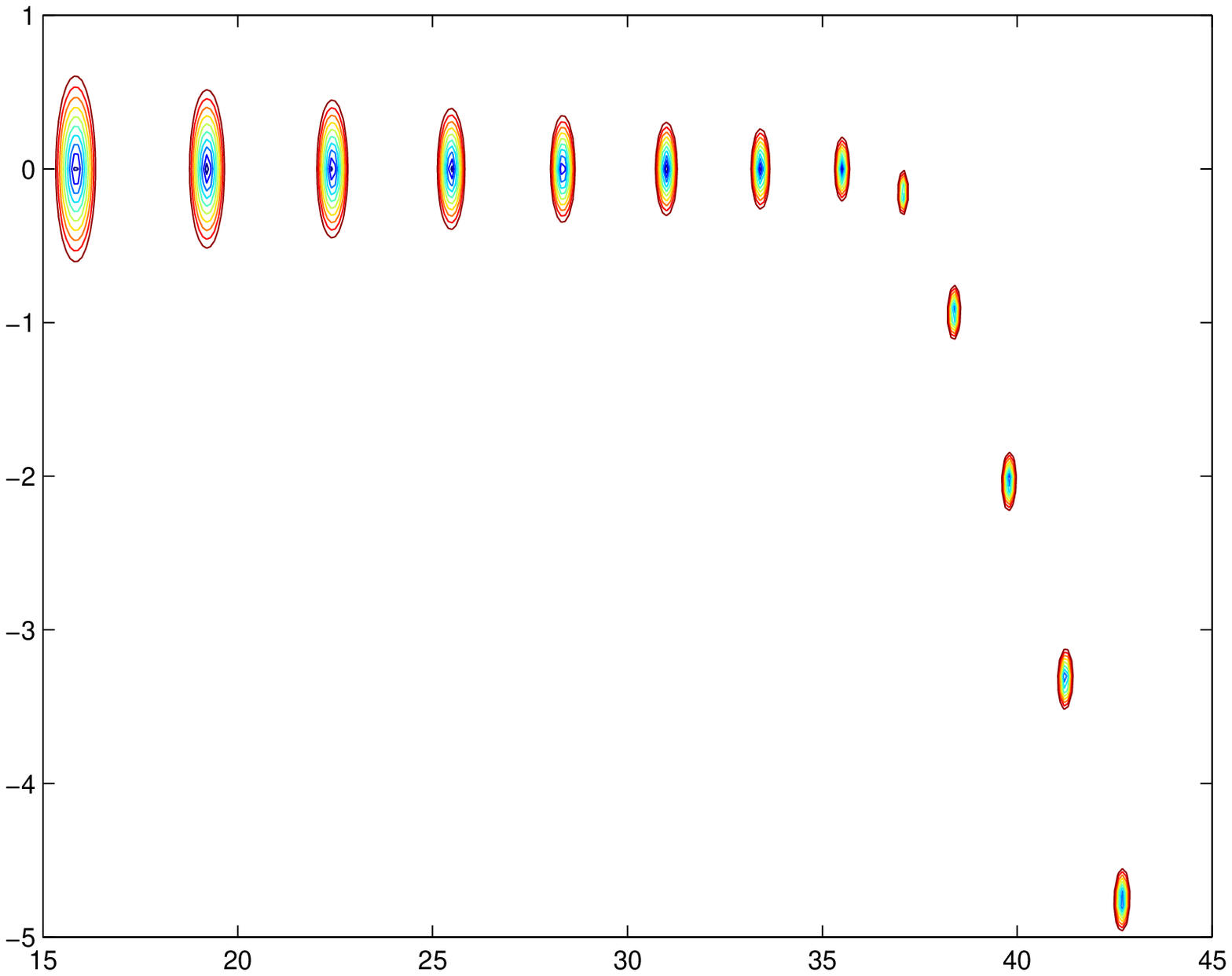}
\end{picture}
\end{figure}

\vspace{2cm}
\begin{center}
Figure 2b. Contour map of $F(z;0)$, $b=10$, $\theta=\pi/4$
\end{center}

\begin{figure}[h]
\begin{picture}(200,200)(0,0)
\includegraphics{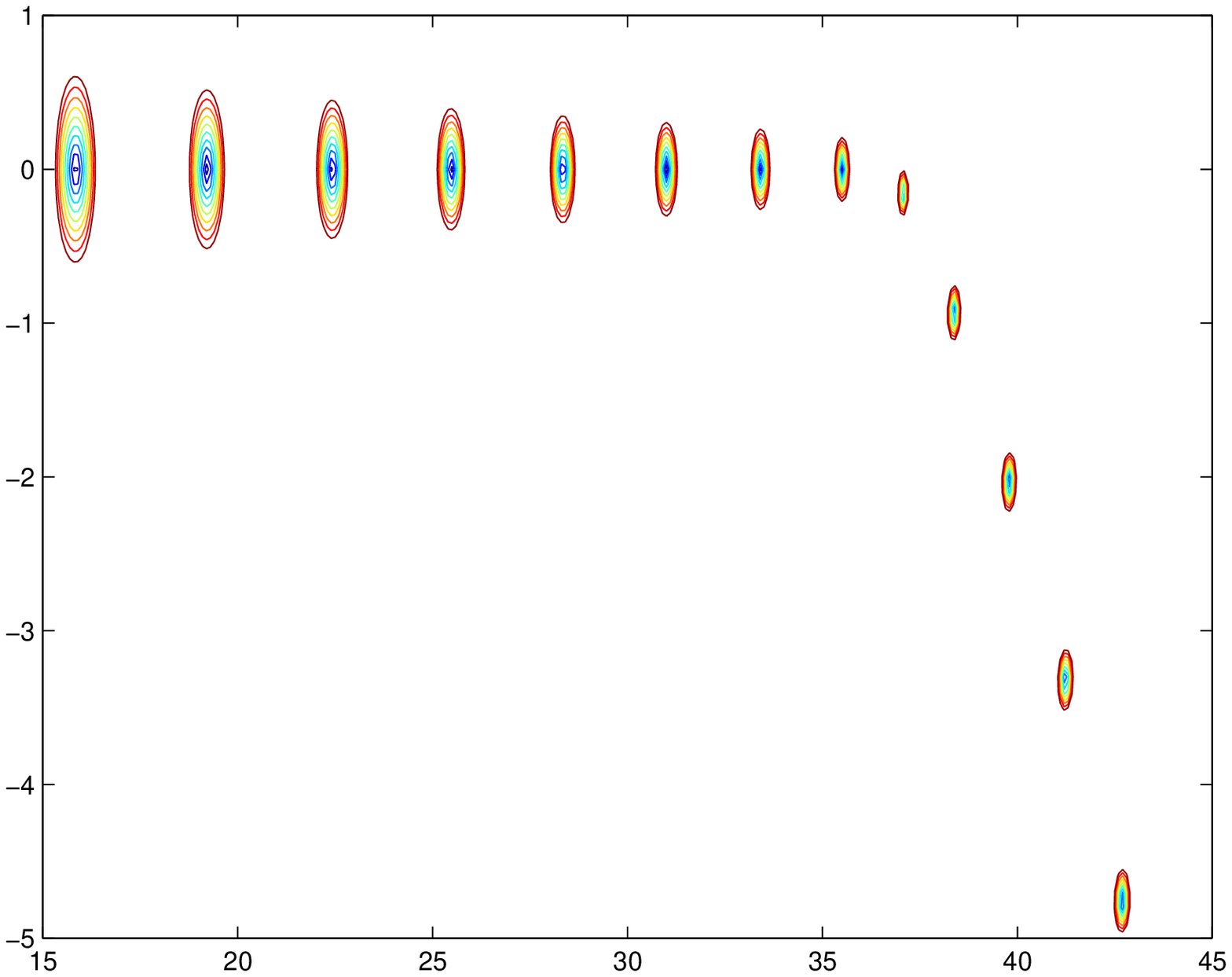}
\end{picture}
\end{figure}

\vspace{1cm}
\begin{center}
Figure 2c. Contour map of $F(z;0)$, $b=10$, $\theta=\pi/16$
\end{center}

\vspace{1cm}


\vskip 1.5in

\begin{thebibliography}{99}
\bibitem{Abr} A. A. Abramov,
Zh. vychisl. Mat. mat. Fiz. {\bf 1}, (1961) 542--545 (Russian).
\bibitem{AbrKon}
A.~A. Abramov, K.~Balla and N.~B.~Konyukhova,
Comput. Math. Banach Center Publs {\bf 13}, (1984) 319--351.
\bibitem{AbrYu}
Abramov A.A. and Yukhno L.F.
Comp. Maths Math. Phys. {\bf 34}, (1994) 671--677.
\bibitem{Bak} N. S. Bakhvalov, {\em Numerical Methods},
(Nauka, Moscow 1973) (Russian).
\bibitem{B1} A. B\"ottcher,
J. Int. Eqns Appl. {\bf 6}, (1994) 267--301.
\bibitem{B2} A. B\"ottcher,
{\em Lectures on Operator Theory and its
Applications} (Fields Institute Monographs, ed. Peter Lancaster.
Amer. Math. Soc. Publ., Providence, RI 1995) 2--74.
\bibitem{CFKS} H. L. Cycon, R. G. Froese, W. Kirsch and B. Simon,
{\em Schr\"odinger Operators; With Application to Quantum Mechanics and
Global Geometry} (Texts and Monographs in Physics,
Springer--Verlag, Berlin 1987).
\bibitem{D1} E. B. Davies {\em Pseudospectra of differential
operators} (Preprint, King's College London, UK 1997).
\bibitem{D2} E. B. Davies,
Proc. Roy. Soc. London A {\bf 454}, (1998) to appear.
\bibitem{D3} E. B. Davies,
Commun. Math. Phys., to appear.
\bibitem{God} S. K. Godunov,
Uspekhi Mat. Nauk {\bf 16}, 3 (99) (1961) 171--174 (Russian).
\bibitem{GK} I. C. Gohberg and M. G. Krein,
Transl. Amer. Math. Soc. {\bf 18}, (1969) 309--316.
\bibitem{HS} P. D. Hislop and I. M. Sigal {\em Introduction to Spectral
Theory} (Springer--Verlag, New York 1966).
\bibitem{Ka}
T. Kato, {\em Perturbation Theory of Linear Operators}
(Springer, Berlin 1966).
\bibitem{Re} S. C. Reddy,
J. Int. Eqns Appl. {\bf 5}, (1993) 369--403.
\bibitem{RedT} S. C. Reddy and L. N. Trefethen,
SIAM J. Appl. Math. {\bf 54}, (1994) 1634--1649.
\bibitem{ReiT} L. Reichel and L. N. Trefethen,
Linear Alg. Appl. {\bf 162--164}, (1992) 153--185.
\bibitem{TT} K.-C. Toh and L. N. Trefethen,
SIAM J. Sci. Comp. {\bf 17}, (1996) 1--15.
\bibitem{T1} L. N. Trefethen,
{\em D F Griffiths and G A Watson, Numerical Analysis 1991}
(Longman Sci. Tech. Publ., Harlow, UK 1992) 234--266.
\bibitem{T2} L. N. Trefethen,
SIAM Review {\bf 39}, (1997) 383--406.
\end{thebibliography}
\end{document}